\documentclass[11pt]{article}
\usepackage{bbm}
\usepackage{mathrsfs}
\usepackage{amscd}
\usepackage{amsmath,amsfonts,amssymb,amscd}
\usepackage{indentfirst,graphics,epsfig,psfrag}
\input{epsf}
\usepackage{ifpdf}
\usepackage{enumerate}
\usepackage{appendix}
\usepackage{enumerate}
\usepackage{color}
\usepackage{lineno}

\makeatletter \@addtoreset{figure}{section} \makeatother
\makeatletter
\long\def\@makecaption#1#2{%
   \vskip 10\p@
   \setbox\@tempboxa\hbox{{#1}\ \ #2}%
   \ifdim \wd\@tempboxa >\hsize

       {#1}\ \ #2\par
   \else
       \hbox to\hsize{\hfil\box\@tempboxa\hfil}%
   \fi}
\makeatother

\newtheorem{thm}{Theorem}[section]
\newtheorem{cor}{Corollary}[section]
\newtheorem{lem}{Lemma}[section]

\newtheorem{obs}{Observation}[section]
\newtheorem{pro}{Proposition}[section]

\newcommand{\qed}{{\hfill\rule{1pt}{1pt}}}

\def\qed{\hfill \rule{4pt}{7pt}}

\begin{document}
\title{\textbf{Steiner $3$-diameter, maximum degree and size\\ of a graph}
\footnote{Supported by the National Science Foundation of China
(Nos. 11551001, 11161037, 11101232, 11461054) and the Science Found
of Qinghai Province (No. 2014-ZJ-907).}}
\author{
\small Yaping Mao\footnote{E-mail: maoyaping@ymail.com}\\[0.3cm]
\small School of Mathematics and Statistics, Qinghai Normal\\
\small University, Xining, Qinghai 810008, China\\[0.2cm]
 }
\date{}
\maketitle
\begin{abstract}
The Steiner $k$-diameter $sdiam_k(G)$ of a graph $G$, introduced by
Chartrand, Oellermann, Tian and Zou in 1989, is a natural
generalization of the concept of classical diameter. When $k=2$,
$sdiam_2(G)=diam(G)$ is the classical diameter. The problem of
determining the minimum size of a graph of order $n$ whose diameter
is at most $d$ and whose maximum is $\ell$ was first introduced by
Erd\"{o}s and R\'{e}nyi. In this paper, we generalize the above
problem for Steiner $k$-diameter, and study the problem of
determining the minimum size of a graph of order $n$ whose Steiner
$3$-diameter is at most $d$ and whose
maximum degree is at most $\ell$.\\[2mm]
\end{abstract}

\section{Introduction}

All graphs in this paper are undirected, finite and simple. We refer
to \cite{Bondy} for graph theoretical notation and terminology not
described here. For a graph $G$, let $V(G)$, $E(G)$ and $e(G)$
denote the set of vertices, the set of edges and the size of $G$,
respectively. We divide our introduction into the following four
subsections to state the motivations and our results of this paper.

\subsection{Distance and its generalizations}

Distance is one of the most basic concepts of graph-theoretic
subjects. If $G$ is a connected graph and $u,v\in V(G)$, then the
\emph{distance} $d(u,v)$ between $u$ and $v$ is the length of a
shortest path connecting $u$ and $v$. If $v$ is a vertex of a
connected graph $G$, then the \emph{eccentricity} $e(v)$ of $v$ is
defined by $e(v)=\max\{d(u,v)\,|\,u\in V(G)\}$. Furthermore, the
\emph{radius} $rad(G)$ and \emph{diameter} $diam(G)$ of $G$ are
defined by $rad(G)=\min\{e(v)\,|\,v\in V(G)\}$ and $diam(G)=\max
\{e(v)\,|\,v\in V(G)\}$. These last two concepts are related by the
inequalities $rad(G)\leq diam(G) \leq 2 rad(G)$. The \emph{center}
$C(G)$ of a connected graph $G$ is the subgraph induced by the
vertices $u$ of $G$ with $e(u)=rad(G)$. Recently, Goddard and
Oellermann gave a survey paper on this subject, see \cite{Goddard}.

The distance between two vertices $u$ and $v$ in a connected graph
$G$ also equals the minimum size of a connected subgraph of $G$
containing both $u$ and $v$. This observation suggests a
generalization of distance. The Steiner distance of a graph,
introduced by Chartrand, Oellermann, Tian and Zou in 1989, is a
natural and nice generalization of the concept of classical graph
distance. For a graph $G(V,E)$ and a set $S\subseteq V(G)$ of at
least two vertices, \emph{an $S$-Steiner tree} or \emph{a Steiner
tree connecting $S$} (or simply, \emph{an $S$-tree}) is a such
subgraph $T(V',E')$ of $G$ that is a tree with $S\subseteq V'$. Let
$G$ be a connected graph of order at least $2$ and let $S$ be a
nonempty set of vertices of $G$. Then the \emph{Steiner distance}
$d_G(S)$ among the vertices of $S$ (or simply the distance of $S$)
is the minimum size among all connected subgraphs whose vertex sets
contain $S$. Note that if $H$ is a connected subgraph of $G$ such
that $S\subseteq V(H)$ and $|E(H)|=d_G(S)$, then $H$ is a tree.
Observe that $d_G(S)=\min\{e(T)\,|\,S\subseteq V(T)\}$, where $T$ is
subtree of $G$. Furthermore, if $S=\{u,v\}$, then $d_G(S)=d(u,v)$ is
nothing new but the classical distance between $u$ and $v$. Set
$d_G(S)=\infty$ when there is no $S$-Steiner tree in $G$.

The distance between two vertices $u$ and $v$ in a connected graph
$G$ also equals the minimum size of a connected subgraph of $G$
containing both $u$ and $v$. This observation suggests a
generalization of distance. The Steiner distance of a graph,
introduced by Chartrand, Oellermann, Tian and Zou in 1989, is a
natural and nice generalization of the concept of classical graph
distance. For a graph $G(V,E)$ and a set $S\subseteq V(G)$ of at
least two vertices, \emph{an $S$-Steiner tree} or \emph{a Steiner
tree connecting $S$} (or simply, \emph{an $S$-tree}) is a such
subgraph $T(V',E')$ of $G$ that is a tree with $S\subseteq V'$. Let
$G$ be a connected graph of order at least $2$ and let $S$ be a
nonempty set of vertices of $G$. Then the \emph{Steiner distance}
$d(S)$ among the vertices of $S$ (or simply the distance of $S$) is
the minimum size among all connected subgraphs whose vertex sets
contain $S$. Note that if $H$ is a connected subgraph of $G$ such
that $S\subseteq V(H)$ and $|E(H)|=d(S)$, then $H$ is a tree.
Observe that $d_G(S)=\min\{e(T)\,|\,S\subseteq V(T)\}$, where $T$ is
subtree of $G$. Furthermore, if $S=\{u,v\}$, then $d_G(S)=d(u,v)$ is
nothing new but the classical distance between $u$ and $v$. Set
$d_G(S)=\infty$ when there is no $S$-Steiner tree in $G$.
\begin{obs}\label{obs1-1}
Let $G$ be a graph of order $n$ and $k$ be an integer with $2\leq
k\leq n$. If $S\subseteq V(G)$ and $|S|=k$, then $d_G(S)\geq k-1$.
\end{obs}

Let $n$ and $k$ be two integers with $2\leq k\leq n$. The
\emph{Steiner $k$-eccentricity $e_k(v)$} of a vertex $v$ of $G$ is
defined by $e_k(v)=\max \{d(S)\,|\,S\subseteq V(G), |S|=k,~and~v\in
S \}$. The \emph{Steiner $k$-radius} of $G$ is $srad_k(G)=\min \{
e_k(v)\,|\,v\in V(G)\}$, while the \emph{Steiner $k$-diameter} of
$G$ is $sdiam_k(G)=\max \{e_k(v)\,|\,v\in V(G)\}$. Note for every
connected graph $G$ that $e_2(v)=e(v)$ for all vertices $v$ of $G$
and that $srad_2(G)=rad(G)$ and $sdiam_2(G)=diam(G)$.

Let $G$ be a $k$-connected graph and $u$, $v$ be a pair of vertices
of $G$. Let $P_k(u,v)=\{P_1,P_2,\cdots,P_k\}$ be a family of $k$
internally vertex-disjoint paths between $u$ and $v$ and
$l(P_k(u,v))$ be the length of the longest path in $P_k(u,v)$. Then
the \emph{$k$-distance} $d_k(u,v)$ between vertices $u$ and $v$ is
the smallest $l(P_k(u,v))$ among all $P_k(u,v)$'s and the
\emph{$k$-diameter} $d_k(G)$ of $G$ is the maximum $k$-distance
$d_k(u,v)$ over all pairs $u,v$ of vertices of $G$. The concept of
$k$-diameter has its origin in the analysis of routings in networks
as described by Chung \cite{Chung}, Du, Lyuu and Hsu \cite{Du}, Hsu
\cite{Hsu, Hsu2}, Meyer and Pradhan \cite{Meyer}.

\subsection{Application backgrounds}

Perhaps the most famous Steiner type problem is the Steiner tree
problem. The original Steiner tree problem was stated for the
Euclidean plane: Given a set of points on the plane, the goal is to
connect these points, and possibly additional points, by line
segments between some pairs of these points such that the total
length of these line segments is minimized. The graph theoretical
version~\cite{Hakimi,Levi} is as follows: Given a graph and a set of
vertices $S$, find a connected subgraph with minimum number of edges
that contains $S$. This is, in general, an $NP$-hard problem
\cite{HwangRW}. There is also a corresponding weighted version.
Obviously, this has applications in computer science and  electrical
engineering. For example, a graph can be a computer network with
vertices being computers and edges being links between them. Here
the Steiner tree problem is to find a subnetwork containing these
computers with the least number of links. We can replace processors
by electrical stations for applications in electrical networks.

Li et al. \cite{LiMaoGutman} gave such a concept. They defined the
{\it $k$-center Steiner Wiener index\/} $SW_k(G)$ of the graph $G$
to be
$$
SW_k(G)=\sum_{\overset{S\subseteq V(G)}{|S|=k}} d(S)\,.
$$
For $k=2$, it coincides with the ordinary Wiener index. One usually
considers $SW_k$ for $2 \leq k \leq n-1$. However, the above
definition can be extended to $k=1$ and $k=n$ as well where
$SW_1(G)=0$ and $SW_n(G)=n-1$. There are other related concepts such
as the Steiner Harary index.  Both indices have chemical
applications~\cite{FurtulaGutmanKatanic, GutmanFurtulaLi}. In
addition, the Steiner degree distance by Gutman \cite{GutmanSDD},
Steiner Harary index by Furtula, Gutman, Katani\'{c}
\cite{FurtulaGutmanKatanic}, Steiner Gutman index by Mao and Das
\cite{MaoDas}, Steiner hyper-Wiener index by Tratnik \cite{Tratnik}
was introduced and studied. We refer the readers to
\cite{FurtulaGutmanKatanic, GutmanFurtulaLi, GutmanSDD, LiMaoGutman,
LiMaoGutman, MaoDas, MaoWangGutmanKlobucar, MaoWangGutmanLi} for
details.

\subsection{Recent progress of Steiner distance}

In \cite{DankelmannSO2}, Dankelmann, Swart and Oellermann obtained a
bound on $sdiam_k(G)$ for a graph $G$ in terms of the order of $G$
and the minimum degree of $G$, that is, $sdiam_k(G)\leq
\frac{3p}{\delta+1}+3n$. Later, Ali, Dankelmann, Mukwembi
\cite{AliDM} improved the bound of $sdiam_k(G)$ and showed that
$sdiam_k(G)\leq \frac{3p}{\delta+1}+2n-5$ for all connected graphs
$G$. Moreover, they constructed graphs to show that the bounds are
asymptotically best possible.

Arunandhi, Cheng and Melekian \cite{ArunandhiChengMelekian} studied
the Steiner $k$-diameters of the tensor product of complete graphs.

The following observation is easily seen.
\begin{obs}\label{obs1-2}
Let $k$ be an integer with $2\leq k\leq n$.

$(1)$ If $H$ is a spanning subgraph of $G$, then $sdiam_k(G)\leq
sdiam_k(H)$.

$(2)$ For a connected graph $G$, $sdiam_k(G)\leq sdiam_{k+1}(G)$.
\end{obs}

In \cite{Chartrand}, Chartrand, Okamoto, Zhang obtained the
following results.
\begin{thm}{\upshape \cite{Chartrand}}\label{th1-1}
Let $G$ be a connected graph of order $n$. Then
$$
k-1\leq sdiam_k(G)\leq n-1.
$$
Moreover, the bounds are sharp.
\end{thm}

\subsection{Classical extremal problem and our generalization}

What is the minimal size of a graph of order $n$ and diameter $d$?
What is the maximal size of a graph of order $n$ and diameter $d$?
It is not surprising that these questions can be answered without
the slightest effort (see \cite{Bollobas}) just as the similar
questions concerning the connectivity or the chromatic number of a
graph. The class of maximal graphs of order $n$ and diameter $d$ is
easy to describe and reduce every question concerning maximal graphs
to a not necessarily easy question about binomial coefficient, as in
\cite{HO1, HS5, O6, W11}. Therefore, the authors study the minimal
size of a graph of order $n$ and under various additional
conditions.

Erd\"{o}s and R\'{e}nyi \cite{ER4} introduced the following problem.
Let $d,\ell$ and $n$ be natural numbers, $d<n$ and $\ell<n$. Denote
by $\mathscr{H}(n,\ell,d)$ the set of all graphs of order $n$ with
maximum degree $\ell$ and diameter at most $d$. Put
$$
e(n,\ell,d)=\min\{e(G):G\in \mathscr{H}(n,\ell,d)\}.
$$
If $\mathscr{H}(n,\ell,d)$ is empty, then, following the usual
convention, we shall write $e(n,\ell,d)=\infty$. For more details on
this problem, we refer to \cite{Bollobas, B26, ER4, ERS1}.

We now consider the generalization of the above problem. Let
$d,\ell$ and $n$ be natural numbers, $d<n$ and $\ell<n$. Denote by
$\mathscr{H}_k(n,\ell,d)$ the set of all graphs of order $n$ with
maximum degree $\ell$ and $sdiam_k(G)\leq d$. Put
$$
e_k(n,\ell,d)=\min\{e(G):G\in \mathscr{H}_k(n,\ell,d)\}.
$$
If $\mathscr{H}_k(n,\ell,d)$ is empty, then, following the usual
convention, we shall write $e_k(n,\ell,d)=\infty$. From Theorem
\ref{th1-1}, we have $k-1\leq d\leq n-1$.

In this paper, we focus our attention on the case $k=3$, and study
the exact value of $e_3(n,\ell,d)$ for $d=n-1,n-2,n-3,2,3$. For
general $d$, we give an upper bound of $e_3(n,\ell,d)$.

\section{Preliminaries}

The following observation is immediate.
\begin{obs}{\upshape \cite{Mao}} \label{obs2-1}
$(1)$ For a cycle $C_n$, $sdiam_k(C_n)=\left
\lfloor\frac{n(k-1)}{k}\right\rfloor$;

$(2)$ For a complete graph $K_n$, $sdiam_k(K_n)=k-1$.
\end{obs}

The following result can be easily proved, which will be used later.
\begin{thm}\label{th2-1}
For $2\leq \ell \leq n-1$ and $3\leq k \leq n$,
$$
e_k(n,\ell,n-1)=n-1.
$$
\end{thm}
\begin{pf}
For $\ell=2$, let $G$ be a path of order $n$. Since $sdiam_k(G)\leq
n-1$, $\Delta(G)=2$ and $e(G)=n-1$, it follows that
$e_k(n,\ell,n-1)\leq n-1$. On the other hand, since we only consider
connected graphs, it follows that $e(G)\geq n-1$ for a connected
graph $G$ of order $n$. So $e_k(n,\ell,n-1)=n-1$, as desired.

Suppose $3\leq \ell \leq n-1$. Let $G$ be a graph obtained from a
star $S_{\ell}$ and a path $P_{n-\ell+1}$ by identifying the center
of the star and one endpoint of the path. Since $\Delta(G)=\ell$,
$sdiam_k(G)\leq n-1$ and $e(G)=n-1$, it follows that
$e_k(n,\ell,n-1)\leq n-1$. On the other hand, since we only consider
the connected graph, it follows that $e(G)\geq n-1$ for a connected
graph $G$ is  of order $n$. So $e_k(n,\ell,n-1)=n-1$. The proof is
now complete.
\end{pf}

\begin{lem}\label{lem2-1}
Let $T$ be a tree of order $n$ with $r$ leaves in $T$. If $r\geq 4$,
then $sdiam_3(T)\leq n-r+2$.
\end{lem}
\begin{pf}
Let $v_1,v_2,\cdots,v_{r}$ be all the leaves of $T$. For any
$S\subseteq V(T)$ and $|S|=3$, there are at least $r-3$ leaves in
$T$ not belonging to $S$. Pick up $r-3$ of these vertices of degree
$1$ and then delete them. The resulting graph is also a tree of
order $n-(r-3)=n-r+3$, say $T'$. By our choosing, it is clear that
$S\subseteq V(T')$, that is, the tree $T'$ is an $S$-Steiner tree.
Therefore, $d_G(S)\leq e(T')=n-r+3-1=n-r+2$. From the arbitrariness
of $S$, we have $sdiam_3(T)\leq n-r+2$. The proof is now complete.
\qed
\end{pf}

\begin{lem}\label{lem2-2}
Let $n,d,\ell$ be three integers with $2\leq d\leq n-2$ and
$n-d+2\leq \ell \leq n-2$. Then
$$
e_3(n,\ell,d)=n-1.
$$
\end{lem}
\begin{pf}
Let $G$ be a graph obtained from a star $S_{\ell}$ and a path
$P_{n-\ell+1}$ by identifying the center of the star and one end of
the path. Clearly, $\Delta(G)=\ell$ and $G$ has $\ell$ leaves. Note
that $\ell\geq n-d+2\geq 4$. From Lemma \ref{lem2-1},
$sdiam_3(G)\leq n-\ell+2\leq n-(n-d+2)+2=d$. Therefore, this graph
shows that $e_3(n,\ell,d)\leq n-1$. On the other hand, we only
consider connected graphs, which implies $e_3(n,\ell,d)\geq n-1$. So
$e_3(n,\ell,d)=n-1$. \qed
\end{pf}\vskip2mm

The following result is from \cite{Mao}.
\begin{lem}{\upshape \cite{Mao}}\label{lem2-3}
Let $G$ be a connected graph of order $n \ (n\geq 3)$. Then
$sdiam_3(G)=n-1$ if and only if $G$ satisfies one of the following
conditions.

$(i)$ $G=T_{a,b,c}$, where $T_{a,b,c}$ is a tree with a vertex $v$
of degree $3$ such that $T_{a,b,c}-v=P_{a}\cup P_{b}\cup P_{c}$,
where $0\leq a\leq b\leq c$.

$(ii)$ $G=C_3(a,b,c)$, where $C_3(a,b,c)$ is a graph containing a
triangle $K_{3}$ and satisfying $C_3(a,b,c)-V(K_{3})=P_{a}\cup
P_{b}\cup P_{c}$, where $0\leq a\leq b\leq c$.
\end{lem}

\begin{cor}\label{cor2-1}
Let $G$ be a connected graph of order $n$. If $sdiam_3(G)\leq n-2$
and $\Delta(G)=2$, then $G$ is cycle of order $n$.
\end{cor}

\begin{lem}{\upshape \cite{Mao}}\label{lem2-4}
Let $G$ be a connected graph of order $n$. Then $sdiam_3(G)=2$ if
and only if $0\leq \Delta(\overline{G})\leq 1$ if and only if
$n-2\leq \delta(G)\leq n-1$.
\end{lem}

\section{For small $d$}

From Theorem \ref{th1-1}, we have $2\leq d\leq n-1$. In this
section, we discuss the cases $d=2$ and $d=3$.

\subsection{The case $d=2$}

If $sdiam_3(G)=2$, then it follows from Lemma \ref{lem2-4} that
$$
n-2\leq \delta(G)\leq n-1, \eqno (4.1)
$$
and hence $n-2\leq \Delta(G)\leq n-1$. So we assume that $n-2\leq
\ell \leq n-1$.
\begin{thm}\label{th3-1}
$(1)$ For $\ell=n-1$, $e_3(n,\ell,2)={n\choose 2}-\frac{n-1}{2}$ for
$n$ odd; $e_3(n,\ell,2)={n\choose 2}-\frac{n-2}{2}$ for $n$ even.

$(2)$ For $\ell=n-2$, $e_3(n,\ell,2)={n\choose 2}-\frac{n}{2}$ for
$n$ even; $e_3(n,\ell,2)=\infty$ for $n$ odd.
\end{thm}
\begin{pf}
$(1)$ For $n$ odd, we let $G$ be a graph obtained from a complete
graph $K_n$ by deleting a maximum matching $M$. Clearly,
$|M|=\frac{n-1}{2}$ and $\Delta(G)=n-1$. From Lemma \ref{lem2-4}, we
have $sdiam_3(G)=2$, and hence $e_3(n,n-1,2)\leq {n\choose
2}-\frac{n-1}{2}$ for $n$ odd. We claim that $e_3(n,n-1,2)=
{n\choose 2}-\frac{n-1}{2}$. Assume, to the contrary, that
$e_3(n,n-1,2)\leq {n\choose 2}-\frac{n-1}{2}-1$. Then there exists a
graph $G$ such that $sdiam_3(G)=2$, $\Delta(G)=n-1$ and $e(G)\leq
{n\choose 2}-\frac{n-1}{2}-1$. Clearly, $\Delta(\overline{G})\geq 2$
and hence $\delta(G)=n-1-\Delta(\overline{G})\leq n-3$, which
contradicts to $(4.1)$. So $e_3(n,n-1,2)={n\choose
2}-\frac{n-1}{2}$.

For $n$ even, let $G$ be a graph obtained from a complete graph
$K_n$ by deleting a matching $M$ such that $|M|=\frac{n-2}{2}$.
Obviously, $\Delta(G)=n-1$. From Lemma \ref{lem2-4}, we have
$sdiam_3(G)=2$, and hence $e_3(n,n-2,2)\leq {n\choose
2}-\frac{n-2}{2}$ for $n$ even. We claim that $e_3(n,n-2,2)=
{n\choose 2}-\frac{n-2}{2}$. Assume, to the contrary, that
$e_3(n,n-1,2)\leq {n\choose 2}-\frac{n-2}{2}-1$. Then there exists a
graph $G$ such that $sdiam_3(G)=2$, $\Delta(G)=n-1$ and $e(G)\leq
{n\choose 2}-\frac{n-2}{2}-1$. Clearly, $\Delta(\overline{G})\geq
1$, $\delta(G)=n-1-\Delta(\overline{G})\leq n-2$, and hence
$\delta(G)=n-2$ by $(4.1)$. Since $e(G)\leq {n\choose
2}-\frac{n}{2}$ and $n$ is even, it follows that $G$ is a graph
obtained from a complete graph $K_n$ by deleting a perfect matching,
which implies $\Delta(G)=n-2$, a contradiction. So
$e_3(n,n-1,2)={n\choose 2}-\frac{n-1}{2}$.

$(2)$ For $n$ even, let $G$ be a graph obtained from a complete
graph $K_n$ by deleting a perfect matching $M$. Obviously,
$|M|=\frac{n}{2}$, $\Delta(G)=n-2$ and $sdiam_3(G)=2$. Therefore,
$e_3(n,n-2,2)\leq {n\choose 2}-\frac{n}{2}$ for $n$ even. We claim
that $e_3(n,n-2,2)= {n\choose 2}-\frac{n}{2}$. Assume, to the
contrary, that $e_3(n,n-2,2)\leq {n\choose 2}-\frac{n}{2}-1$. Then
there exists a graph $G$ such that $sdiam_3(G)=2$, $\Delta(G)=n-2$
and $e(G)\leq {n\choose 2}-\frac{n}{2}-1$. Clearly,
$\Delta(\overline{G})\geq 2$,
$\delta(G)=n-1-\Delta(\overline{G})\leq n-3$, a contradiction. So
$e_3(n,n-1,2)={n\choose 2}-\frac{n}{2}$.

For $n$ odd, let $G$ be a graph such that $sdiam_3(G)=2$ and
$\Delta(G)=n-2$. Since $sdiam_3(G)=2$, it follows that $n-2\leq
\delta(G)\leq n-1$. Since $\Delta(G)=n-2$, it follows that
$\delta(G)=\Delta(G)=n-2$, and hence $G$ is $(n-2)$-regular graph of
order $n$, which is impossible. So $e_3(n,n-2,2)=\infty$ for $n$
odd.\qed
\end{pf}

\subsection{The case $d=3$}

The following proposition and lemma are immediate.
\begin{pro}\label{pro3-1}
Let $K_{n_1,n_2,\cdots,n_r}$ be a complete $r$-partite graph with
$n_1\leq n_2\leq \cdots \leq n_r$. Then
$$
{\rm sdiam}_k(K_{n_1,n_2,\cdots,n_r})=\left\{
\begin{array}{ll}
k-1,&\mbox{{\rm if}~$k>n_r$}\\
k,&\mbox{{\rm if}~$k\leq n_r$.}
\end{array}
\right.
$$
\end{pro}
\begin{pf}
Set $G=K_{n_1,n_2,\cdots,n_r}$. Let $V_1,V_2,\cdots,V_r$ be the
parts of complete $r$-partite graph $G$, and set $|V_i|=n_i \ (1\leq
i\leq r)$. Suppose $k>n_r$. Since $n_1\leq n_2\leq \cdots \leq n_r$,
it follows that for any $S\subseteq V(G)$ and $|S|=k$, there exist
two parts $V_i,V_j$ such that $S\cap V_i\neq \varnothing$ and $S\cap
V_j\neq \varnothing$. Set $S=\{v_1,v_2,\cdots,v_k\}$. Without loss
of generality, let $S\cap V_i=\{v_1,v_2,\cdots,v_s\}$ and $S\cap
V_j=\{v_{s+1},v_{s+2},\cdots,v_{t}\}$. Then the tree induced by the
edges in
$$
\{v_1v_i\,|\,s+1\leq i\leq k\}\cup \{v_{s+1}v_i\,|\,2\leq i\leq s\}
$$
is an $S$-Steiner tree, and hence $d_G(S)\leq k-1$, and hence ${\rm
sdiam}_k(G)\leq k-1$. From Theorem \ref{th1-1}, we have ${\rm
sdiam}_k(G)=k-1$.

Suppose $k\leq n_r$. Choose $S\subseteq V(G)$ and $|S|=k$ such that
$S\subseteq V_r$. Observe that any $S$-Steiner tree must use at
least $k$ edges. Therefore, $d_G(S)\geq k$, and hence ${\rm
sdiam}_k(G)\geq k$. One can easily check that ${\rm sdiam}_k(G)\leq
k$. So ${\rm sdiam}_k(G)=k$, as desired.
\end{pf}

\begin{lem}\label{lem3-1}
Let $T$ be a tree of order $n \ (n\geq 5)$. Then $sdiam_3(T)=3$ if
and only if $T$ is a star.
\end{lem}
\begin{pf}
If $T$ is a star, then $sdiam_3(T)=3$. Conversely, we suppose
$sdiam_3(T)=3$. Suppose that $T$ contains a path $P$ as its subgraph
such that $|V(P)|\geq 3$ and each vertex in $P$ is an internal
vertex of $T$. Let $u,v$ be two endpoints of $P$. Then there exit
two leaves, say $x,y$, such that $xu\in E(T)$ and $yw\in E(T)$.
Choose $S=\{x,u,y\}$. Then $d_G(S)\geq 4$ and hence $sdiam_3(T)\geq
4$, a contradiction. Suppose that $T$ has exactly two internal
vertices of $T$. Since $n\geq 5$, it follows that $T$ has at least
three leaves. Choose three of them as $S$. Then $d_G(S)\geq 4$ and
hence $sdiam_3(T)\geq 4$, a contradiction. We conclude that $T$ has
exactly one internal vertex and hence $T$ is a star.\qed
\end{pf}\vskip1mm

We first give an upper bound of our parameter for general $\ell$ and
$d=3$.
\begin{lem}\label{lem3-2}
For $\frac{n}{2}\leq \ell\leq n-4$,
$$
n\leq e_3(n,\ell,3)\leq
\ell(n-\ell).
$$
\end{lem}
\begin{pf}
Let $K_{\ell,n-\ell}$ be a complete bipartite graph. Since $\ell\geq
\frac{n}{2}$, it follows that $\ell\geq n-\ell$, and hence
$\Delta(G)=\ell$. From Proposition \ref{pro3-1}, we have
$sdiam_3(G)=3$. So $e_3(n,\ell,3)\leq \ell(n-\ell)$. Let $G$ be a
graph with $sdiam_3(G)\leq 3$ and $\frac{n}{2}\leq
\Delta(G)=\ell\leq n-4$. If $sdiam_3(G)=2$, then it follows from
Lemma \ref{lem2-4} that $n-3\leq \delta(G)\leq n-2$, which
contradicts $\Delta(G)=\ell\leq n-4$. Suppose $sdiam_3(G)=3$. From
Lemma \ref{lem3-1}, if $G$ is tree, then $G$ is a star, and hence
$\Delta(G)=n-1$, a contradiction. So $e_3(n,\ell,3)\geq n$.\qed
\end{pf}\vskip1mm

\begin{lem}\label{lem3-3}
For $\ell=n-1$, $e_3(n,\ell,3)=n-1$.
\end{lem}
\begin{pf}
For $\ell=n-1$, let $G$ be a star of order $n$. Then $sdiam_3(G)=3$
and $\Delta(G)=n-1$, and hence $e_3(n,n-1,3)\leq n-1$. Since we only
consider the connected graph, it follows that $e_3(n,n-1,3)\geq
n-1$. So $e_3(n,n-1,3)=n-1$.
\end{pf}\vskip1mm

\begin{lem}\label{lem3-4}
For $\ell=n-2$, $e_3(n,\ell,3)=2n-5$.
\end{lem}
\begin{pf}
For $\ell=n-2$, let $G=K^-_{2,n-2}$ be a graph obtained from a
complete bipartite graph $K_{2,n-2}$ by deleting an edge; see Figure
\ref{figure3-1} $(a)$. Let $u,v,x_{n-2}$ be the vertices of degree
$n-1,n-2,1$ in $K^-_{2,n-2}$, respectively. Set $X=V(G)-
\{u,v,x_{n-2}\}=\{x_1,x_2,\cdots,x_{n-3}\}$.

Now, we show $sdiam_3(G)\leq 3$. It suffices to prove that
$d_G(S)\leq 3$ for any $S\subseteq V(G)$ and $|S|=3$. If $|S\cap
X|=3$, without loss of generality, let $S=\{x_1,x_2,x_3\}$, then the
tree $T$ induced by the edges in $\{ux_1, ux_2,ux_3\}$ is an
$S$-Steiner tree and hence $d_G(S)\leq 3$. If $|S\cap X|=0$, then
$S=\{u,v,x_{n-2}\}$, then the tree $T$ induced by the edges in
$\{ux_1,vx_1,ux_{n-2}\}$ is an $S$-Steiner tree and hence
$d_G(S)\leq 3$. Suppose $|S\cap X|=2$. Then $|S\cap
\{u,v,x_{n-2}\}|=1$. Without loss of generality, let
$S=\{x_1,x_2,u\}$ or $S=\{x_1,x_2,v\}$ or $S=\{x_1,x_2,x_{n-2}\}$.
If $S=\{x_1,x_2,u\}$, then the tree $T_1$ induced by the edges in
$\{ux_1,ux_2\}$ is a Steiner tree connecting $\{x_1,x_2,u\}$. If
$S=\{x_1,x_2,v\}$, then the tree $T_2$ induced by the edges in
$\{vx_1,vx_2\}$ is a Steiner tree connecting $\{x_1,x_2,v\}$. If
$\{x_1,x_2,x_{n-2}\}$, then the tree $T_3$ induced by the edges in
$\{ux_1,ux_2,ux_{n-2}\}$ is a Steiner tree connecting
$\{x_1,x_2,x_{n-2}\}$. Therefore, $d_G(S)\leq 3$. Suppose $|S\cap
X|=1$. Then $|S\cap \{u,v,x_{n-2}\}|=2$. Without loss of generality,
let $S=\{x_1,u,v\}$ or $S=\{x_1,u,x_{n-2}\}$ or
$S=\{x_1,v,x_{n-2}\}$. If $S=\{x_1,u,v\}$, then the tree $T_1$
induced by the edges in $\{ux_1,vx_1\}$ is a Steiner tree connecting
$\{x_1,u,v\}$. If $\{x_1,u,x_{n-2}\}$, then the tree $T_2$ induced
by the edges in $\{vx_1,ux_1, ux_{n-2}\}$ is a Steiner tree
connecting $\{x_1,u,x_{n-2}\}$. If $\{x_1,v,x_{n-2}\}$, then the
tree $T_2$ induced by the edges in $\{ux_1,ux_{n-2}\}$ is a Steiner
tree connecting $\{x_1,v,x_{n-2}\}$. Therefore, $d_G(S)\leq 3$. From
the above argument, we conclude that $d_G(S)\leq 3$ for any
$S\subseteq V(G)$ and $|S|=3$. Since $sdiam_3(G)\leq 3$ and
$\Delta(G)=n-2$, we have $e_3(n,n-2,3)\leq 2n-5$.
\begin{figure}[!hbpt]
\begin{center}
\includegraphics[scale=0.65]{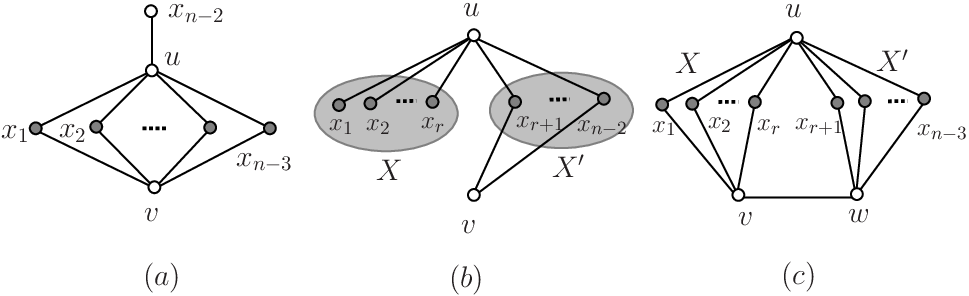}
\end{center}
\caption{Graphs for Lemma \ref{lem3-4}.}\label{figure3-1}
\end{figure}

We now show $e_3(n,n-2,3)\geq 2n-5$. Let $G$ be a graph with
$\Delta(G)=n-2$ and $sdiam(G)\leq 3$. Since $\Delta(G)=\ell=n-2$,
there exists a vertex $u$ such that $d_G(u)=n-2$. Then
$|N_G(u)|=n-2$ and $\{u\}\cup N_G(u)\subseteq V(G)$. Clearly, there
is a vertex $v\in V(G)$ such that $uv\notin E(G)$. Let
$X'=V(G)\setminus \{u,v\}=\{x_1,x_2,\cdots,x_{n-2}\}$. Clearly,
$|X'|\geq 3$. Since $G$ is connected, it follows that there exists a
vertex, say $x_j\in X$ such that $vx_j\in E(G)$. Without loss of
generality, let $vx_{n-2}\in E(G)$. Observe that for each $x_i \
(1\leq i\leq n-3)$ we have $vx_i\in E(G)$ or $vx_i\notin E(G)$.
Without loss of generality, let $vx_1,vx_2,\cdots,vx_r\notin E(G)$
and $vx_{r+1},vx_{r+2},\cdots,vx_{n-3}\in E(G)$, where $0\leq r\leq
n-3$. Set $X=\{x_1,x_2,\cdots,x_r\}$. For each two vertices
$x_i,x_j\in X$, we choose $S=\{x_i,x_j,v\}$. Since $sdiam_3(G)\leq
3$, it follows that any $S$-Steiner tree contains at most $3$ edges.
Therefore, there exists a vertex $x_k\in
\{x_{r+1},x_{r+2},\cdots,x_{n-3}\}$ such that $x_ix_k,x_jx_k\in
E(G)$ or $x_ix_j,x_jx_k\in E(G)$ or $x_ix_j,x_ix_k\in E(G)$. For $r$
even, $e(G)\geq (n-2)+(n-2-r)+r=2n-4$. For $r$ odd, $e(G)\geq
(n-2)+(n-2-r)+(r-1)=2n-5$. So $e_3(n,n-2,3)\geq 2n-5$.

From the above argument, we conclude that $e_3(n,n-2,3)=2n-5$.\qed
\end{pf}\vskip1mm

\begin{lem}\label{lem3-5}
For $\ell=n-3$, $e_3(n,n-3,3)=2n-5$.
\end{lem}
\begin{pf}
Let $G$ be a graph defined as follows; see Figure \ref{figure3-1}
$(c)$.
\begin{eqnarray*}
V(G)&=&\{u,v,w\}\cup \{x_i\,|\,1\leq i\leq n-3\}\\
E(G)&=&\{ux_i\,|\,1\leq i\leq n-3\}\cup \{vx_i\,|\,1\leq i\leq r\}\\
&&\cup \{wx_i\,|\,r+1\leq i\leq n-3\}\cup \{vw\}.\\
\end{eqnarray*}\vskip-0.8cm
\noindent One can check that $sdiam_3(G)=3$ and $\Delta(G)=n-3$.
Therefore, $e_3(n,n-3,3)\leq 2n-5$.

We only need to show $e_3(n,n-3,3)\geq 2n-5$. Let $G$ be a connected
graph with $\Delta(G)=n-3$ and $sdiam_3(G)=3$. It suffices to show
that $e(G)\geq 2n-5$. Since $\Delta(G)=\ell=n-3$, there exists a
vertex $u$ such that $d_G(u)=n-3$, and hence $|N_G(u)|=n-3$ and
$\{u\}\cup N_G(u)\subseteq V(G)$. Clearly, there exist two vertices
$v,w\in V(G)$ such that $v,w\notin N_G(u)$. Let $V(G)\setminus
\{u,v,w\}=\{x_1,x_2,\cdots,x_{n-3}\}$. We have the following two
cases to consider.

\textbf{Case 1.} $vw\in E(G)$

Since $G$ is connected, it follows that there exist two vertices
$x_j,x_k$ such that $vx_j\in E(G)$ and $vx_k\in E(G)$ (note that
$x_j,x_k$ are not necessarily different). Without loss of
generality, let $vx_1,wx_1\in E(G)$ or $vx_1,wx_2\in E(G)$. Suppose
$vx_1,wx_1\in E(G)$. For any $x_i \ (2\leq i\leq n-3)$, we choose
$S=\{x_i,v,w\}$. Since $sdiam_3(G)=3$, it follows that $x_1x_i\in
E(G)$ or $vx_i\in E(G)$ or $wx_i\in E(G)$, and hence $e(G)\geq
(n-3)+2+(n-4)=2n-5$, as desired. Suppose $vx_1,wx_2\in E(G)$. For
$x_3$, we choose $S=\{x_3,v,w\}$. Since $sdiam_3(G)=3$, it follows
that $vx_3,wx_3\in E(G)$ or $x_3x_1,wx_1\in E(G)$ or $x_3x_2,vx_2\in
E(G)$. For any $x_i \ (2\leq i\leq n-3)$, we choose $S=\{x_i,v,w\}$.
Since $sdiam_3(G)=3$, it follows that $x_1x_i\in E(G)$ or $x_2x_i\in
E(G)$ or $x_3x_i\in E(G)$ or $vx_i\in E(G)$ or $wx_i\in E(G)$, and
hence $e(G)\geq (n-3)+4+(n-6)=2n-5$, as desired.

\textbf{Case 2.} $vw\notin E(G)$

Since $G$ is connected, it follows that there exists a vertex $x_j$
such that $vx_j\in E(G)$ or $wx_j\in E(G)$. Without loss of
generality, let $vx_1\in E(G)$. For any $x_i \ (2\leq i\leq n-3)$,
we choose $S=\{x_i,v,w\}$. Since $sdiam_3(G)=3$, it follows that
$x_1x_i\in E(G)$ or $vx_i\in E(G)$ or $wx_i\in E(G)$, and hence
$e(G)\geq (n-3)+1+1+(n-4)=2n-5$, as desired.

From the argument, we conclude that $e(G)\geq 2n-5$, and hence
$e_3(n,n-3,3)\geq 2n-5$.\qed
\end{pf}\vskip1mm

We now conclude our results for $d=3$.
\begin{thm}\label{th3-2}
$(1)$ For $\ell=n-1$, $e_3(n,n-1,3)=n-1$;

$(2)$ For $\ell=n-2$, $e_3(n,n-2,3)=2n-5$;

$(3)$ For $\ell=n-3$, $e_3(n,n-3,3)=2n-5$;

$(4)$ For $\ell=2$, $e_3(n,2,3)=3$ for $n=4$; $e_3(n,2,3)=5$ for
$n=5$; $e_3(n,2,3)=\infty$ for $n\geq 6$.

$(5)$ For $\frac{n}{2}\leq \ell\leq n-4$, $n\leq e_3(n,\ell,3)\leq
\ell(n-\ell)$.
\end{thm}
\begin{pf}
The results in $(1)$-$(3)$ follow from Lemmas \ref{lem3-3},
\ref{lem3-4} and \ref{lem3-5}. Let $G$ be a connected graph with
$\Delta(G)=2$ and $sdiam_3(G)=3$. Then $G=P_n$ or $G=C_n$. If
$G=P_n$, then it follows from Lemma \ref{lem2-3} that
$3=sdiam_3(G)=sdiam_3(P_n)=n-1$, and hence $n=4$. If $G=C_n$, then
it follows from Observation \ref{obs2-1} that
$3=sdiam_3(G)=sdiam_3(C_n)=\lfloor\frac{2n}{3}\rfloor$, and hence
$n=5$. Furthermore, $e_3(n,2,3)=\infty$ for $n\geq 6$. The result in
$(5)$ follow from Lemma \ref{lem3-2}.
\end{pf}

\section{For large $d$}

For $d=n-1$ and $4\leq \ell\leq n-1$, we have proved that
$e_3(n,\ell,n-1)=n-1$. We study the cases $k=n-2,n-3,n-4$ in this
section.

Let $T(a,b,c)$ be a tree obtained from two stars $K_{1,a},K_{1,c}$
and a path $P=v_1v_2\ldots v_b$ by identifying $u$ and $v_1$, $w$
and $v_b$, where $a+b+c=n$, $u$ is the center of $K_{1,a}$, and $w$
is the center of $K_{1,c}$.
\begin{figure}[!hbpt]
\begin{center}
\includegraphics[scale=0.8]{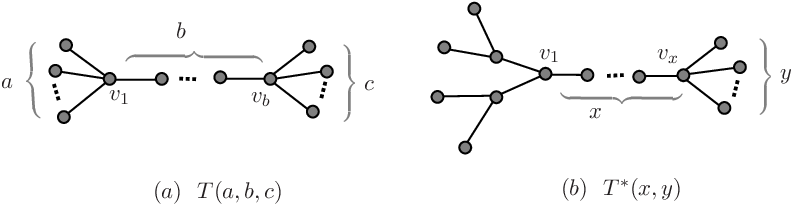}
\end{center}
\caption{Graphs for $n\geq 8$ in Lemma
\ref{lem4-2}.}\label{figure4-1}
\end{figure}
Let $T^*(x,y)$ be a tree obtained from three stars
$K_{1,3},K_{1,3},K_{1,y}$ and a path $P=v_1v_2\ldots v_b$ by
identifying $u,v,v_1$ and $w,v_b$, where $x+y+6=n$, $u$ is a leaf of
$K_{1,3}$, and $v$ is a leaf of another $K_{1,3}$, and $w$ is the
center of $K_{1,c}$.

\subsection{The case $d=n-2$}

In this subsection, we study the case $d=n-2$.
\begin{thm}\label{th4-1}
$(1)$ For $n\geq 4$, $e_3(n,2,n-2)=n$.

$(2)$ For $n\geq 4$,
$$
e_3(n,3,n-2)=\left\{
\begin{array}{ll}
n+1&\mbox{{\rm if}~$n=4$,}\\
n&\mbox{{\rm if}~$n=5$,}\\
n-1&\mbox{{\rm if}~$n\geq 6$.}
\end{array}
\right.
$$

$(3)$ For $n\geq 5$ and $4\leq \ell\leq n-1$, $e_3(n,\ell,n-2)=n-1$.
\end{thm}
\begin{pf}
$(1)$ For $n\geq 4$, let $G=C_n$ be a cycle of order $n$. Then
$\ell=\Delta(G)=2$. From Observation \ref{obs2-1}, we have
$sdiam_3(G)\leq n-2$ and hence this graph shows $e_3(n,2,n-2)\leq
n$. It suffices to show that $e_3(n,2,n-2)\geq n$. Let $G$ be a
graph with $\Delta(G)=2$ and $sdiam_3(G)\leq n-2$. Since
$\Delta(G)=\ell=2$, it follows that $G=P_n$ or $G=C_n$. From Lemma
\ref{lem2-3}, we have $sdiam_3(P_n)=n-1$ and $sdiam_3(C_n)\leq n-2$.
So $e_3(n,2,n-2)\geq n$ and hence $e_3(n,2,n-2)=n$.

$(2)$ For $n\geq 6$, let $G$ be a graph obtained from three paths
$P_1=u_1u_2u_3$, $P_2=w_1w_2w_3$ and $P_3=v_1v_2 \cdots v_{n-4}$ by
identifying the vertices $u_2$ and $v_1$, and then identifying the
vertices $w_2$ and $v_{n-4}$. Clearly, $\ell=\Delta(G)=3$. From
Lemma \ref{lem2-3}, we have $sdiam_3(G)\leq n-2$. Therefore,
$e_3(n,3,n-2)\leq n-1$ and hence $e_3(n,3,n-2)=n-1$. For $n=4$, the
graph $G=K_4^-$ shows that $e_3(4,3,2)\leq 5$. It suffices to show
that $e_3(4,3,2)\geq 5$. Let $G$ be a graph with $\Delta(G)=2$ and
$sdiam_3(G)\leq n-2=2$. Since $\ell=\Delta(G)=3$, there exists a
vertex of degree $3$, say $u_1$. Let $u_1u_2,u_1u_3,u_1u_4\in E(G)$.
Choose $S=\{u_2,u_3,u_4\}$. Since $sdiam_3(G)=2$, it follows that
$u_2u_3,u_3u_4\in E(G)$. Therefore, $e(G)\geq 5$ and hence
$e_3(4,3,2)\geq 5$. So $e_3(4,3,2)=5$. For $n=5$, let $G$ be a graph
obtained from a cycle $C_4$ by adding a pendent edge at one vertex
of $C_4$. One can check that $\Delta(G)=3$ and $sdiam_3(G)\leq 3$.
Therefore, $e_3(5,3,3)\leq 5$. We need to show $e_3(5,3,3)\geq 5$.
Let $G$ be a graph with $\Delta(G)=3$ and $sdiam_3(G)\leq n-2=3$.
Since $\ell=\Delta(G)=3$, there exists a vertex of degree $3$, say
$u_1$. Let $u_1u_2,u_1u_3,u_1u_4\in E(G)$. Since $n=5$, it follows
that there exists a vertex $u_5\in V(G)$. Furthermore, there exists
some vertex $u_j \ (2\leq j\leq 4)$ such that $u_ju_5\in E(G)$.
Without loss of generality, let $u_2u_5\in E(G)$. Choose
$S=\{u_3,u_4,u_5\}$. Since $sdiam_3(G)\leq 3$, it follows that
$u_3u_5\in E(G)$, or $u_4u_5\in E(G)$, or $u_2u_3,u_2u_4\in E(G)$.
Therefore, $e(G)\geq 5$ and hence $e_3(5,3,3)\geq 5$, as desired. So
$e_3(5,3,3)=5$.

$(3)$ For $n\geq 5$ and $4\leq \ell\leq n-1$, let $G$ be a graph
obtained from a star $S_{\ell}$ and a path $P_{n-\ell+1}$ by
identifying the center of the star and one end of the path. Clearly,
$\Delta(G)=\ell$. From Lemma \ref{lem2-3}, we have $sdiam_3(G)\leq
n-2$ and hence this graph shows $e_3(n,\ell,n-2)\leq n-1$. Since we
only consider connected graphs, we have $e_3(n,\ell,n-2)\geq n-1$.
So $e_3(n,\ell,n-2)=n-1$.\qed
\end{pf}

\subsection{The case $d=n-3$}

Let us now turn to the case $d=n-3$.
\begin{lem}\label{lem4-1}
For $n\geq 5$,
$$
e_3(n,2,n-3)=\left\{
\begin{array}{ll}
\infty&\mbox{{\rm if}~$n=5,6$,}\\
n&\mbox{{\rm if}~$n\geq 7$.}
\end{array}
\right.
$$
\end{lem}
\begin{pf}
For $n\geq 7$ and $\ell=2$, let $G=C_n$ be a cycle of order $n$.
From Observation \ref{obs2-1}, we have
$sdiam_3(G)=\lfloor\frac{2n}{3}\rfloor\leq n-3$ and hence this graph
shows $e_3(n,2,n-3)\leq n$. It suffices to prove $e_3(n,2,n-3)\geq
n$. Let $G$ be a connected graph with $\Delta(G)=2$ and
$sdiam_3(G)\leq n-3$. From Corollary \ref{cor2-1}, we have $G=C_n$
and hence $e_3(n,2,n-3)\geq n$. So $e_3(n,2,n-3)=n$. For $n=5,6$, we
have $G=C_n$ by Corollary \ref{cor2-1}. From Observation
\ref{obs2-1}, we have $sdiam_3(G)=\lfloor\frac{2n}{3}\rfloor> n-3$.
So $e_3(n,2,n-3)=\infty$ for $n=5,6$.\qed
\end{pf}

\begin{lem}\label{lem4-2}
For $n\geq 5$,
$$
e_3(n,3,n-3)=\left\{
\begin{array}{ll}
\infty&\mbox{{\rm if}~$n=5$,}\\
n+1&\mbox{{\rm if}~$n=6,7$,}\\
n-1&\mbox{{\rm if}~$n\geq 8$}.
\end{array}
\right.
$$
\end{lem}
\begin{pf}
For $n\geq 8$, let $G=T^*(n-7,1)$. Then $G$ has exactly $5$ leaves
and $\Delta(T')=3$. From Lemma \ref{lem2-1}, $sdiam_3(T')\leq n-3$.
This tree shows that $e_3(n,3,n-3)\leq n-1$ and hence
$e_3(n,3,n-3)=n-1$. For $n=5$, let $G$ be a graph with $\Delta(G)=3$
and $sdiam_3(G)\leq n-3=2$. Furthermore, $sdiam_3(G)=2$. From Lemma
\ref{lem2-4}, we have $3\leq \delta(G)\leq 4$ and hence $3\leq
\delta(G)\leq \Delta(G)=3$, which implies that $G$ is $3$-regular.
The degree sum of graph $G$ is exactly $15$, a contradiction. So
$e_3(5,3,2)=\infty$.

For $n=6$, we let $G=A_4$; see Figure \ref{figure4-1} $(d)$. One can
check that $\Delta(G)=3$ and $sdiam_3(G)\leq n-3=3$. Then
$e_3(6,3,3)\leq 7$. It suffices to show $e_3(6,3,3)\geq 7$. Let $G$
be a graph such that $\Delta(G)=3$ and $sdiam_3(G)\leq 3$. If $G$ is
a tree, then $G=A_1$ or $G=A_2$ or $G=A_3$; see Figure
\ref{figure4-1} $(a),(b),(c)$. Clearly, if we choose vertex set $S$
consisting of three black vertices, then $d_G(S)>3$, which results
in $sdiam_3(G)>3$, a contradiction. So $G$ contains at least one
cycle. Furthermore, we have the following claim.

\textbf{Claim 1.} $G$ contains at least two cycles.

\noindent\textbf{Proof of Claim 1.} Assume, to the contrary, that
$G$ is a unicyclic graph. Let $c(G)$ be the circumference of graph
$G$. Clearly, $3\leq c(G)\leq 6$. If $c(G)=6$, then $G=C_6$, which
contradicts to $\Delta(G)=3$. We may assume that $3\leq c(G)\leq 5$.
If $c(G)=5$, then $G=A_5$; see Figure \ref{figure4-1} $(e)$. One can
also check that $sdiam_3(G)>3$, also a contradiction. If $c(G)=4$,
then $G=A_6$ or $G=A_7$ or $G=A_8$; see Figure \ref{figure4-1}
$(f),(g),(h)$. One can check that $sdiam_3(G)>3$, a contradiction.
For $c(G)=3$, one can also check that $sdiam_3(G)>3$, also a
contradiction.\qed

From Claim 1, $G$ contains at least two cycles and hence $e(G)\geq
7$. So $e_3(6,3,3)\geq 7$, as desired.
\begin{figure}[!hbpt]
\begin{center}
\includegraphics[scale=0.7]{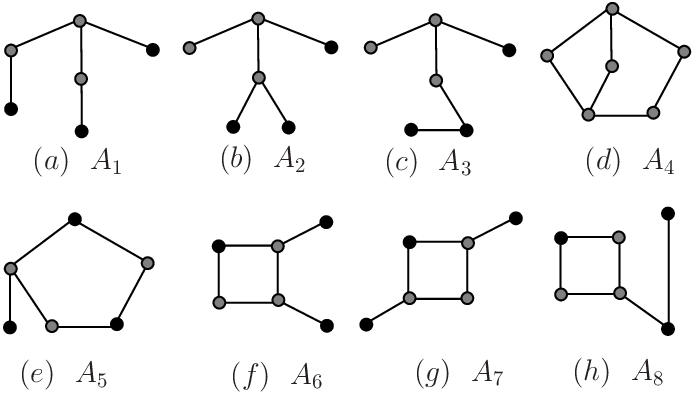}
\end{center}
\caption{Graphs for $n\geq 8$ in Lemma
\ref{lem4-2}.}\label{figure4-1}
\end{figure}

For $n=7$, we let $G=B_1$; see Figure \ref{figure4-2} $(a)$. Then
$\Delta(G)=3$ and $sdiam_3(G)\leq n-3=4$. This graph shows that
$e_3(7,3,4)\leq 8$. It suffices to show that $e_3(7,3,4)\geq 8$. Let
$G$ be connected graph of order $n$ with $sdiam_3(G)\leq 4$ and
$\Delta(G)=3$.

\textbf{Claim 2.} $G$ contains at least two cycles.

\noindent\textbf{Proof of Claim 2.} Assume, to the contrary, that
$G$ is a tree or has exactly one cycle. If $G$ is a tree, then $G\in
\{B_i\,|\,10\leq i\leq 14\}$, and hence $sdiam_3(G)>5$ by choosing
the three black vertices as $S$, a contradiction. We now suppose
that $G$ contains exactly one cycle. Let $c(G)$ denote the
circumference of $G$. Clearly, $3\leq c(G)\leq 7$. If $c(G)=7$, then
$G=C_7$, which contradicts to $\Delta(G)=3$. If $c(G)=6$, then
$G=B_2$; see Figure \ref{figure4-2} $(b)$. By choosing the three
black vertices as $S$, we can see that $sdiam_3(G)>5$, a
contradiction. If $c(G)=5$, then $G=B_3$ or $G=B_4$ or $G=B_5$; see
Figure \ref{figure4-2} $(b)$. By choosing the three black vertices
as $S$, we can see that $sdiam_3(G)>5$, a contradiction. If
$c(G)=4$, then $G\in \{B_i\,|\,6\leq i\leq 9\}$; see Figure
\ref{figure4-2} $(b)$. By choosing the three black vertices as $S$,
we can see that $sdiam_3(G)>5$, a contradiction. For $c(G)=3$, one
can also check that $sdiam_3(G)>5$, also a contradiction.\qed

From Claim 2, $G$ contains at least two cycles, and hence $e(G)\geq
8$. So $e_3(n,3,4)=e(G)\geq 8$, as desired.
\begin{figure}[!hbpt]
\begin{center}
\includegraphics[scale=0.6]{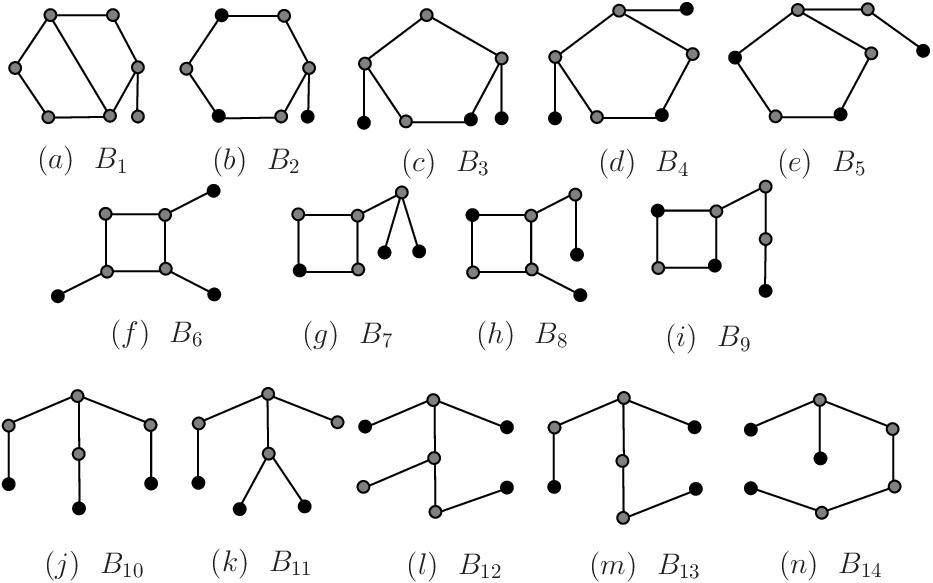}
\end{center}
\caption{Graphs for $n=7$ in Lemma \ref{lem3-2}.}\label{figure4-2}
\end{figure}
\end{pf}

\begin{lem}\label{lem4-3}
For $n\geq 5$,
$$
e_3(n,4,n-3)=\left\{
\begin{array}{ll}
{n \choose 2}-2&\mbox{{\rm if}~$n=5$,}\\
n+1&\mbox{{\rm if}~$n=6$,}\\
n-1&\mbox{{\rm if}~$n\geq 7$}.
\end{array}
\right.
$$
\end{lem}
\begin{pf}
For $n\geq 7$, let $G=T(3,2,2)$. The the tree $T(3,2,2)$ have
exactly $5$ leaves. Clearly, $\Delta(T(3,2,2))=4$. From Lemma
\ref{lem2-1}, $sdiam_3(T(3,2,2))\leq n-3$. This tree shows that
$e_3(n,4,n-3)\leq n-1$, and hence $e_3(n,4,n-3)=n-1$. For $n=5$, let
$G$ be the graph obtained from a complete graph $K_5$ by deleting a
maximum matching. It is clear that $\Delta(G)=4$ and $sdiam_3(G)=2=
n-3$. This graph shows that $e_3(5,4,2)\leq {5 \choose 2}-2=8$. We
need to show that $e_3(5,4,2)\geq 8$. Let $G$ be a graph with
$\Delta(G)=4$ and $sdiam_3(G)\leq n-3=2$. From Lemma \ref{lem2-4},
we have $3\leq \delta(G)\leq 4$. From this together with
$\Delta(G)=4$ and $sdiam_3(G)=2$, it follows that $e(G)\geq 8$ and
hence $e_3(5,4,2)\geq 8$, as desired. So $e_3(5,4,2)=8$.
\begin{figure}[!hbpt]
\begin{center}
\includegraphics[scale=0.7]{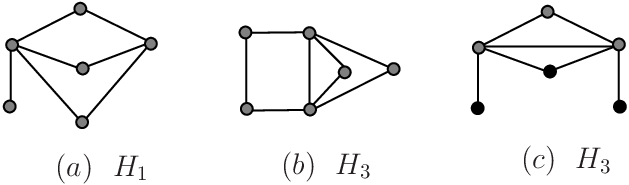}
\end{center}
\caption{Graphs for Lemma \ref{lem4-3}.}\label{figure4}
\end{figure}

For $n=6$, the graph $D_1$ shown in Figure \ref{figure4} $(a)$
satisfies $\Delta(H_1)=4$ and $sdiam_3(D_1)\leq 3=n-3$. This graph
shows that $e_3(6,4,3)\leq 7$. It suffices to prove that
$e_3(6,4,3)\geq 7$. Let $G$ be a graph with $\Delta(G)=4$ and
$sdiam_3(G)\leq 3$. Since $\Delta(G)=4$, it follows that there
exists a vertex of degree $4$, say $u$. Set
$N_G(u)=\{x_1,x_2,x_3,x_4\}$. Since $n=6$, it follows that there is
a remaining vertex in $G$, say $v$. Because $G$ is connected, the
vertex $v$ is adjacent to one of $\{x_1,x_2,x_3,x_4\}$. Without loss
of generality, let $vx_4\in E(G)$. Choose $S=\{x_2,x_3,v\}$. Since
$sdiam_3(G)\leq 3$, it follows that $x_2v\in E(G)$, or $x_3v\in
E(G)$, or $x_2x_4,x_3x_4\in E(G)$. If $x_2v\in E(G)$, then we choose
$S=\{x_1,x_3,v\}$ and hence $x_1v\in E(G)$, or $x_3v\in E(G)$, or
$x_1x_2,x_3x_2\in E(G)$, or $x_1x_4,x_3x_4\in E(G)$. By symmetry, we
only need to consider $x_1v\in E(G)$ or $x_1x_2,x_3x_2\in E(G)$. For
the former, we have $x_1v\in E(G)$ and hence $G=D_1$; see Figure
\ref{figure4} $(b)$. If $x_1x_2,x_3x_2\in E(G)$, then $G=D_2$; see
Figure \ref{figure4} $(c)$. The case $x_3v\in E(G)$ is just the same
as $x_2v\in E(G)$, and so we omit its discussion. We may assume that
$x_2x_4,x_3x_4\in E(G)$. Clearly, $G=D_3$; see Figure \ref{figure4}
$(d)$. Choose the three black vertices to form vertex set $S$. Then
$d_G(S)\geq 4$. So $e_3(6,4,3)=7$ if and only if $G\in \{D_1,D_2\}$.
\qed
\end{pf}\vskip1mm

From Lemmas \ref{lem3-1}, \ref{lem3-2}, \ref{lem3-3} and
\ref{lem2-2}, we have the following result.
\begin{thm}\label{th4-2}
$(1)$ For $n\geq 5$,
$$
e_3(n,2,n-3)=\left\{
\begin{array}{ll}
\infty&\mbox{{\rm if}~$n=5,6$,}\\
n&\mbox{{\rm if}~$n\geq 7$.}
\end{array}
\right.
$$

$(2)$ For $n\geq 5$,
$$
e_3(n,3,n-3)=\left\{
\begin{array}{ll}
\infty,&\mbox{{\rm if}~$n=5$,}\\
n+1&\mbox{{\rm if}~$n=6$,}\\
n&\mbox{{\rm if}~$n=7$,}\\
n-1&\mbox{{\rm if}~$n\geq 8$}.
\end{array}
\right.
$$

$(3)$ For $n\geq 5$,
$$
e_3(n,4,n-3)=\left\{
\begin{array}{ll}
{n \choose 2}-2&\mbox{{\rm if}~$n=5$,}\\
n+1&\mbox{{\rm if}~$n=6$,}\\
n-1&\mbox{{\rm if}~$n\geq 7$}.
\end{array}
\right.
$$

$(4)$ For $n\geq 6$ and $5\leq \ell\leq n-1$, $e_3(n,\ell,n-3)=n-1$.
\end{thm}

\subsection{The case $d=n-4$}

In this subsection, we consider the case $d=n-4$.
\begin{lem}\label{lem4-4}
For $n\geq 5$,
$$
e_3(n,2,n-4)=\left\{
\begin{array}{ll}
\infty&\mbox{{\rm if}~$5\leq
n\leq 9$,}\\
n&\mbox{{\rm if}~$n\geq 10$.}
\end{array}
\right.
$$
\end{lem}
\begin{pf}
For $n\geq 10$ and $\ell=2$, let $G=C_n$ be a cycle of order $n$.
From Observation \ref{obs2-1}, we have
$sdiam_3(G)=\lfloor\frac{2n}{3}\rfloor\leq n-4$ and hence this graph
shows $e_3(n,2,n-4)\leq n$. It suffices to prove $e_3(n,2,n-3)\geq
n$. Let $G$ be a connected graph with $\Delta(G)=2$ and
$sdiam_3(G)\leq n-4$. From Corollary \ref{cor2-1}, we have $G=C_n$
and hence $e_3(n,2,n-4)\geq n$. So $e_3(n,2,n-4)=n$. For $5\leq
n\leq 9$, we have $G=C_n$ by Corollary \ref{cor2-1}. From
Observation \ref{obs2-1}, we have
$sdiam_3(G)=\lfloor\frac{2n}{3}\rfloor> n-4$. So
$e_3(n,2,n-4)=\infty$ for $5\leq n\leq 9$.\qed
\end{pf}

\begin{lem}\label{lem4-5}
For $n\geq 6$,
$$
e_3(n,3,n-4)=\left\{
\begin{array}{ll}
\infty&\mbox{{\rm if}~$n=6$,}\\
n+3&\mbox{{\rm if}~$n=7$,}\\
n+2&\mbox{{\rm if}~$n=8$,}\\
n+1&\mbox{{\rm if}~$n=9$,}\\
n-1&\mbox{{\rm if}~$n\geq 10$}.
\end{array}
\right.
$$
\end{lem}
\begin{pf}
For $n\geq 10$, let $T_3$ be a tree of maximum degree $3$ with
exactly $6$ leaves. Clearly, $\Delta(T'')=3$. From Lemma
\ref{lem2-1}, we have $sdiam_3(T_3)\leq n-4$. This tree shows that
$e_3(n,3,n-4)\leq n-1$ and hence $e_3(n,3,n-4)=n-1$. For $n=6$, let
$G$ be a graph with $\Delta(G)=3$ and $sdiam_3(G)\leq n-4=2$. From
Lemma \ref{lem2-4}, we have $4\leq \delta(G)\leq 5$, which
contradicts to the fact $\Delta(G)=3$. So $e_3(6,3,2)=\infty$.

Suppose $7\leq n\leq 9$. For $n=9$, let $G$ be a graph obtained from
a cycle $C=v_1v_2\ldots v_9$ by adding a new edge $v_1v_5$. One can
easily check that $sdiam_3(G)\leq 5=n-4$ and hence this graph shows
$e_3(9,3,5)\leq 10=n+1$. It suffices to show that $e_3(9,3,5)\geq
10=n+1$. Let $G$ be a graph with $\Delta(G)=3$ and $sdiam_3(G)\leq
5$. Suppose $G$ is a tree. We claim that $G$ has at most $5$ leaves.
Assume, to the contrary, that $G$ has $t \ (6\leq t\leq 8)$ leaves.
Since $\Delta(G)=3$, it follows that $16=2e(G)=\sum_{v\in
V(G)}d(v)\leq t+3(9-t)=27-2t\leq 15$, a contradiction. Since
$\Delta(G)=3$, it follows that $G$ has at least $3$ leaves and at
most $5$ leaves. If $G$ is a tree with $t \ (3\leq t\leq 5)$ leaves,
then we choose three of them as $S$, then $d_G(S)=11-t\geq 6$, which
contradicts to $sdiam_3(G)\leq 5$. We conclude that $G$ is not a
tree, and hence $G$ has cycles. Furthermore, we have the following
claim.

\noindent \textbf{Claim 1.} $G$ contains at least two cycles.

\noindent\textbf{Proof of Claim 1.} Assume, to the contrary, that
$G$ has exactly one cycle. Let $x$ be the number of vertices of
degree $1$ in $G$, and $y$ be the number of vertices of degree $2$
in $G$. We claim that $0\leq x\leq 4$. Assume, to the contrary, that
$x\geq 5$. Since $\Delta(G)=3$, it follows that $18=2e(G)=\sum_{v\in
V(G)}d(v)\leq x+3(9-x)=27-2x\leq 17$, a contradiction. Furthermore,
if $x=4$, then we claim that $y=1$. Assume, to the contrary, that
$y\geq 2$. Then $18=2e(G)=\sum_{v\in V(G)}d(v)\leq
x+2y+3(9-x-y)=27-2x-y=19-y\leq 17$, a contradiction. Similarly, if
$x=3$, then $y\leq 3$.

Suppose $x=4$ and $y=1$. Then $G$ is a unicyclic graph obtained by a
cycle $C=v_1v_2\ldots v_5v_1$ by adding four edges
$v_2v_6,v_3v_7,v_4v_8,v_5v_9$. Choose $S=\{v_6,v_7,v_9\}$. Then
$d_G(S)\geq 6$, which contradicts $sdiam_3(G)\leq 5$. If $x=0$, then
$G=C_9$, and hence $sdiam_3(G)=6$ by Observation \ref{obs2-1}, a
contradiction. If $x=1$, then $G$ is a unicyclic graph obtained by a
cycle $C_r=v_1v_2\ldots v_rv_1$ and a path $P_{9-r}=wv_{r+1}\ldots
v_{9}$ by identifying the vertex $v_1$ in $C_r$ and the endvertex
$w$ in $P_{9-r}$. Choose $S=\{v_{\lfloor
\frac{r}{3}\rfloor},v_{\lfloor \frac{2r}{3}\rfloor},v_9\}$. Since
$3\leq r\leq 8$, it follows that $d_G(S)\geq \lfloor
\frac{2r}{3}\rfloor+(9-r)\geq 6$, which contradicts $sdiam_3(G)\leq
5$.

Suppose $x=3$ and $y\leq 3$. Then $3\leq r\leq 6$. Let $C_r(a,b,c)$
be a graph obtained from a cycle $C_{r}$ and three paths
$P_{a},P_{b},P_{c}$ by adding three edges $z_1u_a,z_2u_b,z_3u_c$,
where $0\leq a\leq b\leq c$, $9=r+a+b+c-3$, $z_1,z_2,z_3$ are three
distinct vertices in $C_{r}$, $u_a,u_b,u_c$ are leaves of
$P_{a},P_{b},P_{c}$, respectively. If $r=3$, then $G=C_3(2,2,2)$ or
$G=C_3(1,2,3)$ or $G=C_3(1,1,4)$. Choose $S$ consisting of all the
three vertices of degree $1$ in $G=C_3(2,2,2)$ or $G=C_3(1,2,3)$ or
$G=C_3(1,1,4)$. Then $d_G(S)\geq 6$, which contradicts
$sdiam_3(G)\leq 5$. If $r=4$, then $G=C_4(1,2,2)$ or $G=C_4(1,1,3)$.
Choose $S$ consisting of all the three vertices of degree $1$ in
$G=C_4(1,2,2)$ or $G=C_4(1,1,3)$. Then $d_G(S)\geq 6$, which
contradicts $sdiam_3(G)\leq 5$. If $r=5$, then $G=C_5(1,1,2)$.
Choose $S$ consisting of all the three vertices of degree $1$ in
$G=C_5(1,1,2)$. Then $d_G(S)\geq 6$, which contradicts
$sdiam_3(G)\leq 5$. If $r=6$, then $G=C_6(1,1,1)$. One can easily
check that $sdiam_3(G)\geq 6$, a contradiction.

Suppose $x=2$. We define two graph classes as follows.
\begin{itemize}
\item Let $\mathcal {G}_9^1$ be a graph class, each graph $G$ of which
is a unicyclic graph obtained by a cycle $C_r=z_1z_2\ldots z_rv_1$
and two paths $P_{s}=u_{1}u_{2}\ldots u_{s},P_{t}=v_{1}v_{2}\ldots
v_{t}$ by identifying a vertex $z_i$ of $C_r$ and the endvertex
$u_{1}$ of $P_{s}$, and then identifying the other vertex $z_j$ of
$C_r$ and the endvertex $v_{1}$ of $P_{t}$, where $9=r+s+t-2$,
$3\leq r\leq 7$ and $1\leq i<j\leq r$.

\item Let $\mathcal {G}_9^2$ be a graph class, each graph $G$ of which
is a unicyclic graph obtained by a cycle $C_r=z_1z_2\ldots z_rz_1$
and $T_{a,b,c}$ (see Lemma \ref{lem2-3}) by identifying a vertex of
$C_r$ and a leaf of $T_{a,b,c}$, where $9=r+a+b+c-3$ and $3\leq
r\leq 6$.
\end{itemize}

For $G\in \mathcal {G}_9^1$, $z_i$ and $z_j$ divide the cycle $C_r$
into two paths $Q_r^1=z_iz_{i+1}\ldots z_j$ and
$Q_r^2=z_{j}z_{j+1}\ldots z_r\ldots z_i$. Without loss of
generality, let $e(Q_r^1)\geq e(Q_r^2)$. Then $e(Q_r^1)\geq
\lceil\frac{r}{2}\rceil$. Choose $S=\{u_{s},v_{t},z_{k}\}$ where
$z_{k}$ is an internal vertex of $Q_r^1$. Since $3\leq r\leq 7$, it
follows that $d_G(S)\geq \lceil\frac{r}{2}\rceil+(9-r)\geq 6$, which
contradicts $sdiam_3(G)\leq 5$.

For $G\in \mathcal {G}_9^2$, let $u,v,w$ be the three leaves in
$T_{a,b,c}$ and $w$ be the identifying vertex in $C_r$. Without loss
of generality, let $w=v_1$. Choose $S=\{u,v,v_{\lfloor
\frac{r}{2}\rfloor}\}$. Then $d_G(S)\geq \lfloor
\frac{r}{2}\rfloor+(9-r)\geq 6$, which contradicts $sdiam_3(G)\leq
5$.\qed

From Claim 1, we conclude that $e_3(9,3,5)=10=n+1$.

For $n=7$, let $G$ be a graph obtained from a cycle $C=v_1v_2\ldots
v_7$ by adding three new edges $v_1v_4,v_2v_5,v_3v_6$. One can check
that $sdiam_3(G)\leq 3=n-4$ and hence this graph shows
$e_3(7,3,3)\leq 10=n+3$. Similarly to the proof of $n=9$, we can
prove that $e_3(7,3,3)\geq 10=n+3$. So $e_3(7,3,3)=10=n+3$. For
$n=8$, let $G$ be a graph obtained from a cycle $C=v_1v_2\ldots v_8$
by adding two new edges $v_1v_5,v_3v_7$. One can check that
$sdiam_3(G)\leq 4=n-4$ and hence this graph shows $e_3(8,3,4)\leq
10=n+2$. Similarly to the proof of $n=9$, we can prove that
$e_3(8,3,4)\geq 10=n+2$. So $e_3(8,3,4)=10=n+2$.
\end{pf}

\begin{lem}\label{lem4-6}
For $n\geq 6$,
$$
e_3(n,4,n-4)=\left\{
\begin{array}{ll}
{n\choose 2}-3&\mbox{{\rm if}~$n=6$,}\\
n+2&\mbox{{\rm if}~$n=7$,}\\
n-1&\mbox{{\rm if}~$n\geq 8$}.
\end{array}
\right.
$$
\end{lem}
\begin{pf}
For $n\geq 8$, let $G=T(3,n-6,3)$. Then $G$ has exactly $6$ leaves
and $\Delta(G)=4$. From Lemma \ref{lem2-1}, $sdiam_3(G)\leq n-4$.
This tree shows that $e_3(n,4,n-4)\leq n-1$ and hence
$e_3(n,4,n-4)=n-1$. For $n=7$, from $(3)$ of Theorem \ref{th3-2}, we
have $e_3(7,4,3)=9$. For $n=6$, let $G$ be a graph with
$\Delta(G)=4$ and $sdiam_3(G)\leq n-4=2$. From Lemma \ref{lem2-4},
we have $4\leq \delta(G)\leq 5$, and hence $G$ is $4$-regular.
Clearly, $G$ is a graph obtained from $K_6$ by deleting a perfect
matching. Therefore, $e(G)={6\choose 2}-3=12$ and
$e_3(6,3,2)=12={n\choose 2}-3$, as desired. \qed
\end{pf}

\begin{lem}\label{lem4-7}
For $n\geq 6$,
$$
e_3(n,5,n-4)=\left\{
\begin{array}{ll}
2n+1&\mbox{{\rm if}~$n=6$,}\\
n+2&\mbox{{\rm if}~$n=7$,}\\
n-1&\mbox{{\rm if}~$n\geq 8$}.
\end{array}
\right.
$$
\end{lem}
\begin{pf}
For $n\geq 8$, let $G=T(4,n-6,2)$. Then $G$ has exactly $6$ leaves
and $\Delta(T_5)=5$. From Lemma \ref{lem2-1}, $sdiam_3(G)\leq n-4$.
This tree shows that $e_3(n,5,n-4)\leq n-1$ and hence
$e_3(n,5,n-4)=n-1$. For $n=7$, from $(2)$ of Theorem \ref{th3-2}, we
have $e_3(7,5,3)=9$. For $n=6$, we let $G$ be the graph obtained
from a complete graph $K_6$ by deleting a matching of size $2$. It
is clear that $\Delta(G)=5$ and $sdiam_3(G)=2=n-4$. This graph shows
that $e_3(6,5,2)\leq {6 \choose 2}-2=13$. We need to show that
$e_3(6,5,2)\geq 13$. Let $G$ be a graph with $\Delta(G)=5$ and
$sdiam_3(G)\leq n-4=2$. From Lemma \ref{lem2-4}, we have $4\leq
\delta(G)\leq 5$. From this together with $\Delta(G)=5$ and
$sdiam_3(G)=2$, it follows that $e(G)\geq 13$ and hence
$e_3(6,5,2)\geq 13$. So $e_3(6,5,2)=13$, as desired.\qed
\end{pf}

\begin{thm}\label{th4-3}
$(1)$ For $n\geq 5$,
$$
e_3(n,2,n-4)=\left\{
\begin{array}{ll}
\infty&\mbox{{\rm if}~$5\leq
n\leq 9$,}\\
n&\mbox{{\rm if}~$n\geq 10$.}
\end{array}
\right.
$$

$(2)$ For $n\geq 6$,
$$
e_3(n,3,n-4)=\left\{
\begin{array}{ll}
\infty&\mbox{{\rm if}~$n=6$,}\\
n+3&\mbox{{\rm if}~$n=7$,}\\
n+2&\mbox{{\rm if}~$n=8$,}\\
n+1&\mbox{{\rm if}~$n=9$,}\\
n-1&\mbox{{\rm if}~$n\geq 10$}.
\end{array}
\right.
$$

$(3)$ For $n\geq 6$,
$$
e_3(n,4,n-4)=\left\{
\begin{array}{ll}
2n&\mbox{{\rm if}~$n=6$,}\\
n+2&\mbox{{\rm if}~$n=7$,}\\
n-1&\mbox{{\rm if}~$n\geq 8$}.
\end{array}
\right.
$$

$(4)$ For $n\geq 6$,
$$
e_3(n,5,n-4)=\left\{
\begin{array}{ll}
2n+1&\mbox{{\rm if}~$n=6$,}\\
n+2&\mbox{{\rm if}~$n=7$,}\\
n-1&\mbox{{\rm if}~$n\geq 8$}.
\end{array}
\right.
$$

$(5)$ For $n\geq 7$ and $6\leq \ell\leq n-1$, $e_3(n,\ell,n-4)=n-1$.
\end{thm}

\section{For general $d$}

We now construct a graph and give an upper bound of $e_3(n,\ell,d)$
for general $\ell$ and $d$.
\begin{pro}\label{th5-2}
For $4\leq d\leq n-1$ and $2\leq \ell\leq n-1$,
$$
e_3(n,\ell,d)\leq \frac{(n-d+1)(n-d+2)}{2}+d-3.
$$
\end{pro}
\begin{pf}
Let $U_p,W_q$ be two cliques of order $p,q$, respectively, where
$p\geq q$ and $p+q=n-d+1$. Set $V(U_p)=\{u_1,u_2,\ldots,u_p\}$ and
$V(W_q)=\{w_1,w_2,\ldots,w_q\}$. Let $G$ be a graph defined as
follow.
\begin{eqnarray*}
V(G)&=&V(U_p)\cup V(W_q)\cup \{v_i\,|\,1\leq i\leq d-1\}\\
E(G)&=&E(U_p)\cup E(W_q)\cup \{v_1u_i:1\leq i\leq p\}\cup \{v_2w_i:1\leq i\leq p\}\\
&&\cup \{u_iw_j:1\leq i\leq p, \ 1\leq j\leq q\}\cup \{v_iv_{i+1}:2\leq i\leq d-2\}\\
\end{eqnarray*}\vskip-0.8cm
It is clear that $|V(G)|=n$, $\Delta(G)=p+q=n-d+1$ and
\begin{eqnarray*}
|E(G)|&=&{p\choose
2}+{q\choose 2}+p+q+pq+d-3\\
&=&\frac{(p+q)^2}{2}+\frac{p+q}{2}+d-3\\
&=&\frac{(n-d+1)(n-d+2)}{2}+d-3.\\
\end{eqnarray*}\vskip-0.8cm
We need to show that $sdiam_3(G)\leq d$. It suffices to prove that
$d_G(S)\leq d$ for any $S\subseteq V(G)$ and $|S|=3$. Set
$X=\{v_1,v_2,\ldots,v_{d-1}\}$. If $S\subseteq X$, the the tree
induced by the edges in $\{v_1u_1,u_1w_1,u_2w_1\}\cup
\{v_iv_{i+1}\,|\,2\leq i\leq d-2\}$ is an $S$-Steiner tree, and
hence $d_G(S)\leq d$. If $|S\cap X|=2$, then $|S\cap V(U_p)|=1$ or
$|S\cap V(W_q)|=1$. Without loss of generality, let $S\cap
V(U_p)=\{u_j\}$. Then the tree induced by the edges in
$\{v_1u_j,u_jw_1,u_2w_1\}\cup \{v_iv_{i+1}\,|\,2\leq i\leq d-2\}$ is
an $S$-Steiner tree, and hence $d_G(S)\leq d$. Suppose $|S\cap
X|=1$. Then $|S\cap V(U_p)|=|S\cap V(W_q)|=1$ or $|S\cap V(U_p)|=2$
or $|S\cap V(W_q)|=2$. We first consider the case $|S\cap
V(U_p)|=|S\cap V(W_q)|=1$. Without loss of generality, let $S\cap
V(U_p)=\{u_j\}$ and $S\cap V(W_q)=\{w_k\}$. If $S=\{u_j,w_k,v_1\}$,
then the tree induced by the edges in $\{u_jv_1,u_jw_k\}$ is an
$S$-Steiner tree, and hence $d_G(S)\leq 2<d$. If $v_1\notin S$, then
the tree induced by the edges in $\{u_jw_k,w_ku_2\}\cup
\{v_iv_{i+1}\,|\,2\leq i\leq d-2\}$ is an $S$-Steiner tree, and
hence $d_G(S)< d$. Next, we consider the case $|S\cap V(U_p)|=2$.
Without loss of generality, let $S\cap V(U_p)=\{u_j,u_k\}$. If
$S=\{u_j,u_k,v_1\}$, then induced by the edges in
$\{u_jv_1,u_kv_1\}$ is an $S$-Steiner tree, and hence $d_G(S)\leq
2<d$. If $v_1\notin S$, then the tree induced by the edges in
$\{u_kw_1,u_jw_1,w_1u_2\}\cup \{v_iv_{i+1}\,|\,2\leq i\leq d-2\}$ is
an $S$-Steiner tree, and hence $d_G(S)\leq d$. In the end, we
consider the case $|S\cap V(W_q)|=2$. Without loss of generality,
let $S\cap V(W_q)=\{w_j,w_k\}$. If $S=\{u_j,u_k,v_1\}$, then induced
by the edges in $\{w_ju_1,w_ku_1,v_1u_1\}$ is an $S$-Steiner tree,
and hence $d_G(S)\leq 3<d$. If $v_1\notin S$, then the tree induced
by the edges in $\{w_kv_2,w_jv_2\}\cup \{v_iv_{i+1}\,|\,2\leq i\leq
d-2\}$ is an $S$-Steiner tree, and hence $d_G(S)<d$. We conclude
that $e_3(n,\ell,d)\leq \frac{(n-d+1)(n-d+2)}{2}+d-3$. \qed
\end{pf}

\end{document}